\documentclass[12pt]{amsart}
\usepackage{mathrsfs}
\usepackage{amssymb}

\usepackage[centertags]{amsmath}
\usepackage{amsfonts}
\usepackage{amsthm}
\usepackage{times}

\usepackage{enumerate} %
\usepackage{xcolor}
\usepackage[margin=1.25in]{geometry}
\usepackage{centernot}
\usepackage{csquotes}
\usepackage[symbol]{footmisc}

\newtheorem{theorem}{Theorem}[section]
\newtheorem*{theorem*}{Theorem}

\newtheorem{corollary}[theorem]{Corollary}
\newtheorem{lemma}[theorem]{Lemma}
\newtheorem{remark}{Remark}
\newtheorem{proposition}[theorem]{Proposition}
 
\newtheorem{problem}{Problem}[section]

\newtheorem{question}{Question}[section]
\theoremstyle{definition}

\newcommand{\nn}{\mathbb{N}}

\newcommand{\rr}{\mathbb{R}}
\newcommand{\ext}{\mathrm{ext}}
\newcommand{\supp}{\mathrm{supp}}

\DeclareMathOperator{\sgn}{sgn}

\newcommand{\Isom}{\mathrm{Isom}}
\newcommand{\Id}{\mathrm{Id}}
\newcommand{\card}{\mathrm{card}}

\begin{document}

	\title{Surjective isometries on Banach sequence spaces: A survey}
	
	\author{Leandro Antunes}
	\author{Kevin Beanland}
	
	\address{Departamento de Matem\'atica, Universidade Tecnol\'ogica Federal do Paran\'a, Campus Toledo \\
		85902-490 Toledo, PR \\ Brazil}
	\email{leandroantunes@utfpr.edu.br}

	\address{Department of Mathematics, Washington and Lee University, Lexington, VA 24450.}
	\email{beanlandk@wlu.edu}

	\thanks{2010 \textit{Mathematics Subject Classification}. Primary: }
	\thanks{\textit{Key words}: }
	
	
	\begin{abstract}
		In this survey, we present several results related to characterizing the surjective isometries on Banach sequence spaces. Our survey includes full proofs of these characterizations for the classical spaces as well as more recent results for combinatorial Banach spaces and Tsirelson-type spaces. Along the way, we pose many open problems related to the structure of the group of surjective isometries and characterizations of the group of surjective isometries for various Banach spaces. 
	\end{abstract}
	
	\maketitle
	
	\tableofcontents
	
	\section{Introduction}
	
	The study of norm-preserving linear operators goes back to the beginning of Operator Theory. According to Fleming and Jamison \cite{fleming-isometries}, the first characterization of isometries on a particular space may have been by I. Schur \cite{schur} in 1925, in connection with matrix spaces. The Banach-Stone Theorem \cite{Ba-book, Stone}  is a milestone in this area.
	
	\begin{theorem}
		Let $K$ and $L$ be compact Hausdorff spaces and $T:C(K)\to C(L)$ be a surjective linear isometry. Then there is a homeomorphism $\phi:K \to L$ and $g\in C(L)$ with $|g(y)|=1$ for all $y\in L$ so that $(Tf)(y)=g(y)f(\phi(y))$ for all $y\in L$ and $f \in C(K)$.
	\end{theorem}
	
	This theorem was originially proved by Banach for compact metric spaces and subsequently extended by Stone to the currently stated form (for a further extension see Sherman in \cite{Sherman-BStone})). Banach, in his famous treatise \cite{Ba-book} in 1932, exhibits the general form of surjective linear isometries $T: X \to X$, which he called \textit{rotations}, of several classical Banach spaces $X$, such as $C[0,1]$ and $L_p[0,1]$ with $p \geqslant 1$.

	In the current paper, we survey many results characterizing the group of surjective linear isometries on various Banach spaces. Our emphasis is on Banach spaces whose underlying vectors space consists of infinite sequences of real numbers.
	
	In his monograph, Banach does not present the details for the proofs of the characterizations of isometries of sequence spaces, such as $c$, $c_0$, and $\ell_p$, stating that they are analogous to the function spaces. One of the purposes of this paper is to exhibit the complete proofs for the characterization of these and many other groups of isometries. A sketch of the proof for the classical spaces $c_0$ and separable $\ell_p$ can be found in Lindenstrauss and Tzafriri \cite[Proposition 2.f.14]{LTz-book}. A characterization of non-surjective isometries is due to Lamperti \cite{Lamperti-Isom} and can be found in the recent monograph \cite{CMR-Book}. In more recent work, Diestel and Koldobsky consider the case of Sobolov spaces \cite{Diestel-Sobolov,DiestelKold-Positivity}.
	
	In the past few decades, several authors have characterized all surjective linear isometries on different classical and semi-classical sequence spaces including Schreier's space, James' space, and Tsirelson's space. There is a vast literature studying linear isometries on Banach spaces including the aforementioned book of Fleming and Jamison. 
	
	In this current survey, we focus our attention on real spaces although in many cases the complex space has corresponding results. The goal of the survey is to present a mixture of classical and new results related to the structure of the group of surjective isometries on a Banach space. 
	
	By a \textit{Banach sequence space} we will mean a Banach space of scalar sequences for which the coordinate functionals are continuous. Let $c_{00}$ be the set of all finitely supported scalar sequences indexed by $\mathbb{N}$. If $x=(a_i) \in c_{00}$ let $\text{supp}\,x= \{i \in \mathbb{N}: a_i\not=0\} $. Let $(e_i)$ be the standard unit vectors in $c_{00}$. For a real Banach space $X$ denote by $\Isom(X)$ the group of all surjective linear isometries $T: X \to X$. All sequence spaces contain $c_{00}$ as a (usually dense) subset and in many cases $(e_i)$ will be a Schauder basis for $X$. A Schauder basis $(x_i)$ for a Banach space $X$ is 1-unconditional if $\|\sum \varepsilon_i a_i x_i\|\leqslant\|\sum  a_i x_i\| $ for $\varepsilon_i \in\{\pm 1\}$ and $(a_i) \in c_{00}$.
	
	A sequence space $X$ has a \textit{standard isometry group} if $\Isom(X)$ consists of all isometries of the form
	$$T(a_n) = (\varepsilon_n a_{\pi(n)})$$
	where $(\varepsilon_n)$ is a sequence in $\{-1,1\}$ and $\pi \in \text{SYM}(\mathbb{N})$, i.e., $\pi: \nn \to \nn$ is a permutation of natural numbers. If $X$ has a standard isometry group it is isomorphic to the wreath product of $\mathbb{Z}_2$ and $\text{SYM}(\mathbb{N})$ denoted $\mathbb{Z}_2 \,\wr\, \text{SYM}(\mathbb{N})$. The elements are pairs $((\varepsilon_i)_i,\pi)$, $((\delta_i)_i,\sigma)$ in $\mathbb{Z}_2^\mathbb{N}\times \text{SYM}(\mathbb{N})$ and the product is given by
	$$((\varepsilon_i)_i,\pi)\cdot ((\delta_i)_i,\sigma)=(\varepsilon_i\delta_{\pi^{-1}(i)},\pi \circ \sigma).$$  
	$X$ has a \textit{diagonal standard isometry group} if the only permutation allowed is the identity; that is, $\Isom(X)$ is isomorphic to $\mathbb{Z}_2^\mathbb{N}$. A Banach space $X$ has a \textit{trivial} isometry group if the only isometries are $\{\pm \Id\}$. 
	
	In this survey, we will give examples of sequences spaces having standard isometry groups, diagonal-standard isometry groups, and trivial isometry groups. As we will see, there are in-between cases that arise naturally. For a Banach space $(X,\|\cdot\|)$ a norm $\|\cdot \|_1$ is an equivalent norm (or a renorming) of $X$ if there are $A,B>0$ so that $A\|x\|\leqslant \|x\|_1\leqslant B\|x\|.$ For a given Banach space $X$ and a group $G$ whether there is an equivalent norm on $X$ whose group of isometries is isomorphic to $G$ is a fascinating question that has been studied by many authors \cite{FGalego-TAMS,FR-Extracta, FR-Duke}. We will preview some of this literature and state some open problems. In general, for the sake of clarity of exposition, we have frequently chosen to present less general forms of the results in the literature. This survey is far from comprehensive. We hope the reader will be persuaded by the presentation to peruse the many citations and work that has been done in this area.

	\section{Isometries on the classical sequence spaces.}
	
	With the notable exception of $\ell_2$, each of the classical spaces $c_0,c$, and $\ell_p$ for $p\in[1,\infty]\setminus\{2\}$ have a standard isometry group. We first present the case of $\ell_2$ which is the same for any separable Hilbert space. 	This well-known result illustrates the divergence between the structure of surjective isometries on $\ell_2$ and the other classical spaces. 
	
	\begin{theorem} Let $T:\ell_2 \to \ell_2$ be a bounded linear operator. Then
		$T \in \Isom(\ell_2)$ if and only for any orthonormal basis $(x_n)$ of $\ell_2$ there is an orthonormal basis $(y_n)$ so that $Tx_n=y_n$. \label{ell2}
	\end{theorem}

	\begin{proof}
		Suppose that $T \in \Isom(\ell_2)$ and let $(x_n)$ be an orthonormal basis of $\ell_2$. Recall the following identity for $x,y\in \ell_2$, 
		\begin{equation}
		2\langle x,y\rangle  = \|x\|^2-\|x-y\|+\|y\|^2. \label{Hilbert}
		\end{equation} 
		It follows that $(Tx_n)$ is an orthonormal basis for $\ell_2$. 
		
		Suppose $T:\ell_2 \to \ell_2$ is a bounded linear map so that $(Tx_n)$ and $(x_n)$ are orthonormal bases. Using (\ref{Hilbert}) for each $n \in \mathbb{N}$, $\langle Tx_n,Tx\rangle= \langle x_n,x\rangle.$ Then from Parseval's theorem
		$\|Tx\|^2= \sum_n |\langle Tx_n,Tx\rangle|^2 = \sum_n |\langle x_n,x\rangle|^2 = \|x\|^2.$  \end{proof}
	
	Let $B_X$ and $S_X$ denote the closed unit ball and unit sphere of a Banach space $X$, respectively. In particular, Theorem \ref{ell2} shows that any Hilbert space $H$ has a \textit{transitive} norm in the sense that for each $x \in S_H$ the orbit of $x$ under the action $\Isom(H)$ is $S_H$. The longstanding \textit{Mazur Rotation Problem} asks whether any separable Banach space with such a norm must be isometric to a separable Hilbert space \cite{Ba-book}. Dilworth and Randrianantoanina \cite{DilRand-JFA} recently proved the striking result that if $X$ is an $\ell_p$ space for $1<p<\infty$ and $p\not=2$, there is no equivalent norm on $X$ so that $X$ is almost transitive in this equivalent norm. That is, there exists $x \in S_X$ so that the orbit of $x$ with respect to the group of isometries under the norm is not dense in $S_X$ (and therefore, trivially, not equal $S_X$). Since $\Isom(\ell_p)$ is a standard isometry group for $1\leqslant p<\infty$ and $p\not=2$ it is easy to see that the standard norm is not almost transitive. The result of Dilworth and Randrianantoanina says this remains true for any equivalent renorming of such an $\ell_p$. Moreover, in this paper, they show that the same conclusion holds for subspaces and quotients of such $\ell_p$ spaces. To contrast this result note that, in 1979, Lusky proved that each (separable) Banach space is a $1$-complemented subspace of a (separable) Banach space having an almost transitive norm \cite{Lusky-Studia}. 
	For an excellent exposition of many results related to Mazur's Rotation Problems see \cite{CFR-Mazur}. 
	
	We now prove that $\ell_p$ for $p \in[1,\infty)\setminus\{2\}$ has a standard isometry group.  We first show that such isometries take disjointly supported vectors to disjointly supported vectors. This proof can be found in Banach's book (see \cite[p. 175] {Ba-book}), where the proofs are presented for $L^p$ spaces.
	
	\begin{lemma}
		\label{lemma24}
		Let $1 \leqslant p < \infty$, $p \neq 2$, and let $x$ and $y \in \ell_p$ such that $\supp \, x \cap \supp \, y = \emptyset$. If $T \in \Isom(\ell_p)$, $z = T(x)$ and $w = T(y)$, then $\supp \, z \cap \supp \, w = \emptyset$. \label{disjoint}
	\end{lemma}
	
	\begin{proof}
		Let $x=(a_i)$, $y=(b_i)$, $x=(c_i)$ and $w=(d_i)$. Since $\supp \,x \cap \supp\, y = \emptyset$, for every $\alpha, \beta \in \rr$ we have
		$$\|\alpha x + \beta y\|^p = |\alpha|^p\|x\|^p+|\beta|^p\|y\|^p$$
		and since $T$ is an isometry,
		\begin{equation}
		\label{eq2}
		\|\alpha z + \beta w\|^p = |\alpha|^p\|z\|^p+|\beta|^p\|w\|^p
		\end{equation}
		We will divide the rest of the proof in three cases:
		
		Suppose $p = 1$: Taking $\alpha = \beta = 1$ in \eqref{eq2}, we have
		$$\sum_{i=1}^\infty |c_i+d_i| = \sum_{i=1}^\infty |c_i|+ \sum_{i=1}^\infty |d_i|.$$
		On the other hand, taking $\alpha = 1$ and $\beta = -1$ in \eqref{eq2}, we have
		$$\sum_{i=1}^\infty |c_i-d_i| = \sum_{i=1}^\infty |c_i|+ \sum_{i=1}^\infty |d_i|.$$
		Notice that if $c_i \neq 0$ and $d_i \neq 0$ then we must have $|c_i+d_i| < |c_i|+|d_i|$ or $|c_i-d_i| < |c_i|+|d_i|$. Hence, the conditions above are only satisfied when $\supp\, z \cap \supp\, w = \emptyset$.
		
		Suppose now $2 < p < \infty$: Let $H = \supp \, z \cap \supp \, w$ and let
		\begin{eqnarray*}
			\varphi: \rr \times \nn & \to & \rr \\
			(\alpha,n) & \mapsto & |\alpha c_n  + \beta d_n|^p.		
		\end{eqnarray*}
		Then
		$$\dfrac{\partial \varphi}{\partial \alpha}(\alpha,n) = p|\alpha c_n  + \beta d_n|^{p-1}c_n \sgn(c_n)$$
		and
		\begin{equation}
		\label{eq5}
		\dfrac{\partial^2 \varphi}{\partial \alpha^2}(\alpha,n) = p(p-1)|\alpha c_n  + \beta d_n|^{p-2}c_n^2.
		\end{equation}
		Since $(|\alpha c_n  + \beta d_n|^{p-1}) \in \ell_{\frac{p}{p-1}}$ and $(c_n) \in \ell_p$, by Hölder's inequality we have $\sum_{n=1}^\infty |\alpha c_n  + \beta d_n|^{p-1}|c_n| < \infty$. Therefore, using \eqref{eq2} we have
		\begin{eqnarray*}
			\sum_{n\in H} \dfrac{\partial \varphi}{\partial \alpha}(\alpha,n) & = & \dfrac{\partial}{\partial \alpha} \sum_{n\in H} \varphi(\alpha,n)\\
			& = & \dfrac{\partial}{\partial \alpha} \sum_{n \in H} |\alpha c_n  + \beta d_n|^p\\
			& = & \dfrac{\partial}{\partial \alpha} \left(|\alpha|^p \sum_{n \in H} |c_n|^p  + |\beta|^p \sum_{n \in H} |d_n|^p\right)\\
			& = & p|\alpha|^{p-1} \sgn(\alpha) \sum_{n \in H} |c_n|^p.\\
		\end{eqnarray*}
		Hence,
		$$\sum_{n \in H} \dfrac{\partial^2 \varphi}{\partial \alpha^2}(\alpha,n) = \dfrac{\partial}{\partial \alpha} \left(\sum_{n\in H} \dfrac{\partial \varphi}{\partial \alpha}(\alpha,n)\right) = p(p-1)|\alpha|^{p-2} \sum_{n \in H} |c_n|^p.$$
		Using \eqref{eq5} we may conclude that
		$$|\alpha|^{p-2} \sum_{n \in H} |c_n|^p = \sum_{n \in H} |\alpha c_n  + \beta d_n|^{p-2}c_n^2.$$
		Taking $\alpha = 0$ and $\beta = 1$,
		$$\sum_{n \in H} |d_n|^{p-2}c_n^2 = 0$$
		and since $|d_n|^{p-2}, c_n^2 \geqslant 0$, we must have $H = \emptyset$.
		
		Case 3 is $1 < p < 2$: We will use a dual argument. Let $q$ be the conjugate of $p$. Again, we claim that $z$ and $w$ have disjoint support. Let $z^*=(|c_n|^{p-1}\sgn(c_n))_n$ and $w^*=(|d_n|^{p-1}\sgn(d_n))_n$ and note that these vectors have the same supports as $z$ and $w$ respectively. By definition,
		$$\langle z^*,z\rangle = \|z^*\|_q\|z\|_{p}.$$
		Since $T$ is an isometry
		$$\langle T^*z^*,x\rangle = \langle z^*,z\rangle = \|z^*\|_q\|z\|_{p}= \|T^*z^*\|_q\|x\|_{p}.$$
		Notice third equality above is equality in H\"older's inequality. This implies that $T^*z^*$ is a multiple of $x$ and, therefore, they have the same support. The same holds for $T^*w^*$ and $y$. Since $T^*$ is an isometry on $\ell_q$ and $q>2$ we apply the previous case to conclude that the supports of $z^*$ and $x^*$ are disjoint which suffices to prove our claim. \end{proof}
	
	\begin{theorem}
		\label{teor73}
		Let $1 \leqslant p < \infty$ and $p \not=2$. Then, $\ell_p$ has a standard group of isometries. 
	\end{theorem}
	
	\begin{proof}
		For each $n \in \mathbb{N}$, let $A_n=\text{supp}\,Te_n$. By Lemma \ref{disjoint}, $A_n\cap A_m=\emptyset$. Since $T$ is surjective $\cup A_n = \mathbb{N}$. The theorem will follow when we show $A_n$ is a singleton for each $n$. Find the sequence $(b_j)_{j=1}^\infty$ so that for each $n$, $Te_n= \sum_{j \in A_n} b_j e_j$ and notice the $b_j\not=0$ for each $j$ since $T $ is surjective. Fix $n_0$ and $j_0\in A_{n_0}$. We claim that $\{j_0\}=A_{n_0}$. Let $x=\sum_n a_n e_n$ so that $Tx=e_{j_0}$. Then 
		$$e_{j_0}= \sum_{n=1}^\infty a_n T e_n = \sum_{n\not=n_0} a_n (\sum_{j \in A_n} b_j e_j ) + a_{n_0} (\sum_{j \in A_{n_0}} b_j e_j).$$
		For $n\not=n_0$ we have $a_nb_j=0$ for $j \in A_n$. Since $b_j\not=0$ for each $j$ we have $a_n=0$ for $n\not=n_0$. For $j \in A_{n_0}$ and $j\not=j_0$ we have $a_{n_0}b_j=0$ since $a_{n_0}\not=0$, there is no such $j\not=j_0$. This is the desired result. 
	\end{proof}
	
	Let $c$ denote the Banach space of all convergent sequences of real numbers equipped with the supremum norm. The space $c$ is isomorphic but not isometric to $c_0$ and the duals of $c$ and $c_0$ are both isometric to $\ell_1$. 
	
	\begin{corollary}
		\label{cor26}
		The spaces $c$ and $c_0$ each have a standard isometry group.
	\end{corollary}
	
	\begin{proof}
		We prove only case of $c$ as the case of $c_0$ is identical. Let $T \in \Isom(c)$, $x = (a_i)_{i=1}^\infty \in c$ and $y = (b_i)_{i=1}^\infty = T(x)$. Since $T \in \Isom(c)$ then $T^* \in \Isom(\ell_1)$. There is $(\varepsilon_i)$ in $\{-1,1\}$ and a permutation $\pi: \nn \to \nn$ such that
		$T^*(c_i) = (\varepsilon_i c_{\pi(i)})_{i=1}^\infty$ for every $z = (c_i)_{i=1}^\infty \in \ell_1$.
		
		Let $(e_i)_{i=1}^\infty$ and $(e_i^*)_{i=1}^\infty$ denote the standard unit vector basis of $c$ and $\ell_1$, respectively. Recall that the elements $e_k^* \in \ell_1$ act on $x \in c$ by the formula $e_k^*(x) = a_k.$
		Hence, for every $k \in \nn$, $$\varepsilon_ka_{\pi(k)}=\varepsilon_ke_{\pi(k)}^*(x)=T^*(e_k^*)(x) = e_k^*(T(x)) = e_k^*(y)=b_k.$$
		This is the desired result
	\end{proof}
	
	Essentially the same proof can be adapted for other $\ell_1$ preduals $X$ so long as the isometry between $X^*$ and $\ell_1$ has the required form. Another proof of this fact for $c_0$, found in \cite[Chapter 5, Theorem 10]{lax}, and that we reproduce below with minor modifications, uses the fact that isomorphisms preserve the codimension of closed subspaces. 
	
	\begin{proposition}
		\label{prop27}
		The isometry group of $c_0$ is standard.
	\end{proposition}
	
	\begin{proof}
		It suffices to check that for every $j \in \nn$, $T(e_j) = \pm e_k$ for some $k \in \nn$. For each $j \in \nn$ let
		$$A_j = \{(a_n) \in c_0 \colon a_j = 0\}.$$
		Notice that $A_j$ is a $1$-codimensional closed subspace of $c_0$. Therefore, $B_j = T(A_j)$ is also a $1$-codimensional closed subspace of $c_0$. 
		
		We claim that $B_j = A_k$ for some $k \in \nn$. Indeed, let $f_j = T(e_j)$. Since $T$ is an isometry, for every $v \in B_j$ such that $\|v\| \leqslant 1$, we have
		\begin{equation}
		\label{eq1}
		\|f_j + v\| = \|e_j + T^{-1}(v)\| = 1,
		\end{equation}
		because $T^{-1}(v) \in A_j$ and $\|T^{-1}(v)\| \leqslant 1$. Moreover, $\|f_j\| = \|e_j\| = 1$ and hence, there exists $k \in \nn$ such that $f_j(k)= \pm 1$. By equation \eqref{eq1} we conclude that $v(k)=0$, for every $v \in B_j$ and, since $B_j$ is a $1$-codimensional space of $c_0$, we must have $B_j = A_k$, which implies that $f_j = \pm e_k$.
	\end{proof}

	\subsection{Extreme Points and Isometries} A vector $x$ in a convex set $S$ is called an extreme point of $S$ if $x$ does not lie in the interior of any nontrivial closed line segment of $X$. We denote by $\ext(S)$ the set of extreme points of $S$. The following equalities are standard exercises: $$\ext(B_{c_0})=\emptyset,\, \ext(B_{\ell_1})=\{\pm e_n: n \in \mathbb{N}\}, \ext(B_{\ell_p}) =S_{\ell_p}\text{ for }p\in(1,\infty)$$ and $$\ext(B_{\ell_\infty}) = \{ (\varepsilon_i)_i: \varepsilon_i \in \{\pm 1\}\}.$$ 
	In many cases, the standard fact that a surjective isometry on a space $X$ takes $\ext(B_X)$ onto $\ext(B_X)$ is exploited to characterize all isometries. Indeed, the reader should observe that this can be used to give an easy proof that $\ell_1$ has a standard group of isometries. 
	The cardinality of $\ext(B_X)$ has been studied in several cases. Lindenstrauss and Phelps \cite{LP-Israel} proved that if $X$
	is reflexive then $\ext(B_X)$ is uncountable. This proof uses a Baire category theory argument and is non-constructive. Consequently, it can be more difficult to characterize $\Isom(X)$ for reflexive $X$ by simply using extreme point arguments. On the other end of the spectrum, in \cite{BDHQ} the authors show that $\ext(B_X)$ is countable for all combinatorial spaces $X$ (see the next section for the definition). This builds on work by Shura and Trautman for Schreier's space \cite{ShTr-Glasgow}.
	
	As an illustration of the usefulness of characterizing the extreme points of a Banach space we prove that $\ell_\infty$ has a standard group of isometries.
	
	\begin{proposition}
		The space $\ell_\infty$  has a standard group of isometries.\label{ellinfty}
	\end{proposition}
	
	\begin{proof}
		Let $U \in \Isom(\ell_\infty)$. For each $A\subset \mathbb{N}$ let $w_A \in \ext(B_{\ell_\infty})$ be the vector that is $1$ on $A$ and $-1$ on $\mathbb{N}\setminus A$. Since $U(w_{\nn})\in \ext(B_{\ell_\infty})$, then $U(w_{\nn})=w_B$ for some $B$. Therefore, if $V_1$ is the isometry sending $w_B \mapsto w_\mathbb{N}$ we have that $V_1U(w_{\nn})=w_{\nn}$. Let $U_1=V_1U$. 
		
		We show there  is a permutation $\pi$ so that $U_1w_A=w_{\pi(A)}$ for each $A$. Define $\tilde{U}:\mathcal{P}(\mathbb{N}) \to \mathcal{P}(\mathbb{N})$ by  $U_1(w_A)=w_{\tilde{U}(A)}$. We claim that $\tilde{U}$ is an order preserving automorphism of $\mathcal{P}(\nn)$. From this it will follow that $\tilde{U}(A)=\pi(A)$ for some permutation $\pi$.
		Indeed, note that for $A,B,C \subset \mathbb{N}$ we have 
		$$w_A+w_B-w_C=w_\nn$$
		holds if and only if $A\cap B=C$. Applying $U$ to both sides we have $w_{\tilde{U}(A)}+w_{\tilde{U}(B)}-w_{\tilde{U}(C)}=w_\mathbb{N}$ and so $\tilde{U}(A\cap B)=\tilde{U}(A)\cap\tilde{U}(B)$. That is, $\tilde{U}$ preserves intersections and hence preserves subset inclusion, as desired.  Let $\pi$ so that $U_1 w_A=w_{\pi(A)}$. Let $V_2$ be the isometry induced by $\pi^{-1}$. Hence, $V_2U_1w_A=w_A$ for all $A \subset \mathbb{N}$. Let $U_2=V_2U_1$. 
		
		We claim $U_2$ is the identity. From this, it will follow that $U$ has the desired form. For each $A\subset \nn$ define $e_{A}=1/2(w_\nn +w_A)$ (the indicator vector). As $U_2$ is linear it fixes all such vectors and hence fixes all sequences with finitely many distinct elements. The collection of sequences with finitely many values is dense in $B_{\ell_\infty}$. Whence $U_2$ is the identity. 
	\end{proof}

	\begin{remark}
		\textup{The above proof was provided by Yves de Cornulier in \cite{YCor}. Essentially the same proof works in the case of $\ell_\infty(X)$ of any set $X$.} 
	\end{remark}
	
	\begin{remark}
		\textup{Theorem \ref{ellinfty} is also a consequence of the Banach-Stone theorem. As $\ell_\infty$ is isometric to $C(\beta\mathbb{N})$ where $\beta\mathbb{N}$ is the Stone-{\v C}ech compactification of $\mathbb{N}$, any surjective isomorphism of $\ell_\infty$ corresponds to a self-homeomorphism of $\mathbb{N}$ which must take isolated points to isolated points. Since such a map must be induced by a permutation; the result follows. }
	\end{remark}

	\section{Isometries on combinatorial spaces}
	
	In this section, we consider another class of sequence spaces, called combinatorial spaces. The most well-known example of a combinatorial space is Schreier's space. Schreier sets arose naturally in Banach space theory to solve a classical problem and have become important tools in numerous constructions of Banach sequence spaces. These sets were independently discovered in combinatorics in relation to coloring problems for finite subsets of $\mathbb{N}$.

	Banach and Saks proved in \cite{Banach-Saks} that every bounded sequence in $L^p$, with $p > 1$, has a subsequence such that its arithmetic means converge in norm. This property is currently called the \textit{Banach-Saks property}. They asked whether this property is also valid for the space of continuous functions. This question was answered negatively by J. Schreier in \cite{Schreier}, by constructing a sequence of functions $(f_i)$ in $C[0,1]$ that converges weakly to $0$ but has no subsequence whose arithmetic means converge in norm. In his construction, Schreier considered finite sets of natural numbers $F$ such that $\card F \leqslant \min F$. The family
	$$\mathcal{S}_1 = \{F \in [\nn]^{<\infty}: \card F \leqslant \min F\} \cup \{\emptyset\}$$
	is now called \textit{Schreier family} and its elements are called the \textit{Schreier sets}.
	
	Nishiura and Waterman showed in \cite{NishiuraWaterman} that every space with the Banach-Saks property is reflexive. On the other hand, Baernstein \cite{Baernstein} gave an example of a reflexive space that does not have the Banach-Saks property. In \cite{Beauzamy} Beauzamy used a variation of Baernstein's construction to define the space $X_{\mathcal{S}_1}$ as the completion of $c_{00}$ in respect to the norm
	$$\|x\|_{X_{\mathcal{S}_1}} = \sup_{F \in \mathcal{S}_1} \sum_{k \in F} |a_k|, \quad x = (a_1,\dots,a_k,\dots) \in c_{00}.$$ This space was called \textit{Schreier space}, \index{Schreier! space} in \cite{BeauzamyLapreste}. Beauzamy proved that the interpolation space $(\ell_1,X_{\mathcal{S}_1})_{\theta,p}$, for $0 < \theta < 1$ and $1 < p < \infty$ is a reflexive space that does not have the Banach-Saks property, providing simpler examples than the one obtained by Baernstein.
	
	If $E,F \subset \mathbb{N}$ we write $E<F$ if $\max E < \min F$ and if $x,y \in c_{00}$ we write $x<y$ if $\text{supp}\, x < \text{supp}\, y$. In their landmark paper \cite{AlA-Dissertationes}, Alspach and Argyros defined \textit{higher order Schreier set} as follows. Let $\mathcal{A}_n$ denote the set of finite subsets of $\nn$ with cardinality less than or equal to $n$. Letting $\mathcal{S}_0=\mathcal{A}_1$ and supposing that $\mathcal{S}_\alpha$ has been defined for some countable ordinal $\alpha$, we define
	$$\mathcal{S}_{\alpha+1}=\{\bigcup^{n}_{i=1} E_i: n \leqslant E_1 < E_2 < \dots < E_n \mbox{ and }E_i \in  \mathcal{S}_\alpha\}\cup \{\emptyset\}.$$
	If $\alpha$ is a limit ordinal then we fix $\alpha_n\nearrow\alpha$ and define $$\mathcal{S}_\alpha=\{\emptyset\} \cup\{F \subset \nn: \text{ for some } n \geqslant 1,  F\in\mathcal{S}_{\alpha_{n}} \text{ and } n \leqslant F\}.$$  We say that $\mathcal{S}_\alpha$ is the \textit{Schreier family of order $\alpha$}. \index{Schreier! family of order $\alpha$} The 
	\textit{Schreier space of order $\alpha$}, denoted by $X_{\mathcal{S}_\alpha}$, is defined as the completion of $c_{00}$ with respect to the norm
	$$\|x\|_{X_{\mathcal{S}_\alpha}} = \sup_{F \in \mathcal{S}_\alpha} \sum_{k \in F} |a_k|, \quad x = (a_1,\dots,a_k,\dots) \in c_{00}.$$
	
	By transfinite induction that each family $\mathcal{S}_\alpha$ is:
	
	\begin{enumerate}
		\item \textit{hereditary} ($F \in \mathcal{S}_\alpha$ and $G \subset F \implies G \in \mathcal{S}_\alpha$);
		\item \textit{spreading} ($\{l_1,l_2,\dots,l_n\} \in \mathcal{S}_\alpha$ and $l_i \leqslant k_i \implies \{k_1,k_2,\dots,k_n\} \in \mathcal{S}_\alpha$);
		\item \textit{compact}, seen as a subset of $\{0,1\}^\nn$, with the natural identification of $\mathcal{P}(\nn)$ with $\{0,1\}^\nn$.
	\end{enumerate}
	
	A collection $\mathcal{F}$ of finite subsets of $\nn$ satisfying these three properties is called a \textit{regular family} in \cite{BDHQ}.  An element $F$ in a regular family $\mathcal{F}$ is \textit{maximal} \index{Maximal set} if there is no set $G \in \mathcal{F}$ that contains $F$ properly. We denote by $\mathcal{F}^{\max}$ the collection of maximal sets of $\mathcal{F}$.
	For a regular family $\mathcal{F}$ that contains all singletons, the \textit{combinatorial Banach space} $X_{\mathcal{F}}$ is defined to be the completion of $c_{00}$ with respect to the norm
	$$\|x\|_{X_{\mathcal{F}}} = \sup_{F \in \mathcal{F}} \sum_{k \in F} |a_k|, \quad x = (a_1,a_2,\dots) \in c_{00}.$$
	If $\mathcal{F}$ is regular family  containing the singletons, the spaces $X_\mathcal{F}$ and $C(\mathcal{F})$ (continuous functions on $\mathcal{F}$) both encode topological information about $\mathcal{F}$. The set $\mathcal{F}$ is homeomorphic to a subset of the dual unit sphere in the weak$^*$ topology. Moreover, the Banach-Stone theorem guarantees the topological structure of $\mathcal{F}$ encodes the isometric structure of $C(\mathcal{F})$.
	
	The term combinatorial Banach space was first defined by Gowers in his blog \cite{Go-blog}, for a hereditary set of finite subsets of natural numbers containing the singletons. There is some inconsistency in the literature regarding the definition of combinatorial family, however, at minimum, a combinatorial family should contain all singletons and be closed under subsets. That is, the spreading and compact assumptions are sometimes excluded. In the case that the family $\mathcal{F}$ is compact, the corresponding space $X_\mathcal{F}$ has a shrinking Schauder basis and is, therefore, more closely aligned with $c_0$ and the classical Schreier space.

	Unlike the classical sequence spaces, the unit vector basis of Schreier's space is not a $1$-symmetric Schauder basis. This is easily seen by considering the permutation switching $1$ and $2$ and fixing all $n>2$. Indeed, the only permutation $\pi$ so that $\|\sum_{i}a_i e_i\|=\|\sum_i a_i e_{\pi(i)}\|$ is the identity. The same also holds for $X_{\mathcal{S}_\alpha}$ for any non-zero countable $\alpha$. Note, however, this is not true for any combinatorial space, as $X_{\mathcal{S}_0}$ is isometric to $c_0$. 
	
	We now give a mostly self-contained proof that Schreier's space $X_{\mathcal{S}_1}$ has a diagonal standard group of isometries. Although Schreier spaces have been extensively studied and used fruitfully in the literature in the past 50 years, the following result appeared only recently \cite{ABC-Illinois}. 
	
	\begin{theorem}
		\label{Isom}
		The space $X_{\mathcal{S}_1}$ has a diagonal standard group of isometries.
	\end{theorem}
	
	\begin{proof}
		Let $T \in \Isom(X_{\mathcal{S}_1 })$. For each $i \in \mathbb{N}$, let
		$$T(e_i) = f_i = (f_i(1),f_i(2), f_i(3), \dots)$$
		and
		$$T^{-1}(e_i) = d_i = (d_i(1), d_i(2), d_i(3), \dots).$$
		We will divide the proof of this theorem in six steps:

		\textit{Step 1: $f_1 = \pm e_1$.}
		
		\noindent Assume towards a contradiction that there exists $k \in \supp f_1$ with $k>1$ and let $\ell > k$. Since $\{k,\ell\} \in\mathcal{S}_1$ we have
		$$\|e_1 \pm d_\ell\| = \|f_1 \pm e_\ell\| \geqslant |f_1(k)| + 1>1.$$
		However, since $\|d_\ell\|=1$, if $F \in \mathcal{S}_1$ satisfies
		$$\sum_{i \in F}|e_1(i) + d_{\ell}(i)|>1$$
		then $1 \in F$ and so $F=\{1\}$. Likewise if $G \in \mathcal{S}_1$ satisfies
		$$\sum_{i \in G}|e_1(i) - d_{\ell}(i)|>1$$
		Then $G=\{1\}$. Consequently
		$$|1 + d_{\ell}(1)|>1 \mbox{ and } |1 - d_{\ell}(1)|>1.$$
		This cannot happen since $|d_{\ell}(1)|\leqslant 1$. Therefore there is no $k\geqslant 2$ with $f_{1}(k)\not=0$. Hence, $f_1=\pm e_1$, as desired.
		
		\textit{Step 2: For every $i \in \mathbb{N}$, $d_i, f_i \in c_{00}$ and there exists a non maximal Schreier set $F_i$ such that $\sum_{k \in F_i}|f_i(k)| = 1$}.
		
		\noindent In \cite{BDHQ}, the authors prove that if $x \in \ext(B_{\mathcal{S}_n})$ then $x \in c_{00}$ and there exists a non-maximal Schreier set $F$ such that $\sum_{k \in F}|x(k)| = 1$. A simple computation show that that $e_1 + e_i \in \ext(B_{X_{\mathcal{S}_1}})$ for every $i \geqslant 2$. Since isometries preserve extreme points, and $f_1=\pm e_1$, we have $f_1 + f_i$, $-f_1 + f_i$, $f_1-f_i$, $d_1 + d_i$, $-d_1 + d_i$, $d_1-d_i \in \ext(B_{\mathcal{S}_1})$, which gives us the desired conclusion. 
		
		\textit{Step 3: For every $i > \max \supp f_2$, $f_i = \pm e_m$ for some $m \geqslant 2$}
		
		\noindent Let $n = \max \supp f_2$ and let $F$ be a non-maximal Schreier set such that $\sum_{k\in F} |f_2(k)| = 1$. Then, $\{n+1\} \cup F \in \mathcal{S}_1$ and 
		$$\|e_2 + d_i\| = \|T(e_2 + d_i)\| = \|f_2 + e_i\| = 2.$$
		Hence, there exists $m \geqslant 2$ such that $d_i(m) = \pm 1$. Since $\|d_i\| = \|e_i\| = 1$, we must have $d_i = \pm e_m$, for some $m \geqslant 2$.
		
		\textit{Step 4: For every $N \in \nn$ there exist $i, j > N$ such that $f_i = \pm e_j$}
		
		\noindent  It follows from Step 3 and the injectivity of $T$.
		
		\textit{Step 5: for every $i > \max \supp f_2$, $f_i = \pm e_i$}
		
		\noindent By Step 3, $f_i = \pm e_j$ for some $j \geqslant 2$. Assume towards a contradiction that $j \neq i$. Without loss of generality, suppose that $j > i$. By Step 4, there exist $\alpha_1, \alpha_2, \dots, \alpha_i, \beta_1, \beta_2, \dots, \beta_i > j$ such that $f_{\alpha_n} = \pm e_{\beta_n}$, for every $n \in \{1, \dots, i\}$. Then,
		$$T(e_i + \sum_{n=1}^i e_{\alpha_n}) = e_j + \sum_{n=1}^i \pm e_{\beta_n}.$$
		On the other hand, $\|e_i + \sum_{n=1}^i e_{\alpha_n}\| = i$ and $\|e_j + \sum_{n=1}^i \pm e_{\beta_n}\| = i+1$, which contradicts the fact that $T$ is an isometry.
		
		\textit{Step 6: $f_2 = \pm e_2$}
		
		\noindent Suppose that $\max \supp f_2 = k > 2$. Then, by Step 5 we have
		$$\|e_2 + \sum_{n=k+1}^{2k+1} e_n\| = \|f_2 + \sum_{n=k+1}^{2k+1} \pm e_n\|.$$
		However, $\|e_2 + \sum_{n=k+1}^{2k+1} e_n\| = k$, while $\|f_2 + \sum_{n=k+1}^{2k+1} \pm e_n\| \geqslant |f_2(k)| + k > k$. Therefore, we must have $\max \supp f_2 \leqslant 2$. Moreover, if $f_2(1) \neq 0$, then $\|e_1 + f_2\| > 1$ or $\|e_1 - f_2\| > 1$, but $\|T^{-1}(e_1)\pm T^{-1}(f_2)\| = \|\pm e_1 \pm e_2\| = 1$, which contradicts the fact that $T^{-1}$ is an isometry. Hence, we must have $f_2(1)=0$ and, since $\|f_2\|=1$, it follows that $f_2(2) = \pm 1$.	
	\end{proof}
	
	In \cite{ABC-Illinois}, the above method of proof is generalized to show that $X_{\mathcal{S}_n}$ for each $n \geqslant \mathbb{N}$, has a diagonal standard isometry group.  Shortly after this theorem appeared, Brech, Ferenczi, and Tcaciuc \cite{BFT-PAMS} proved for each non-zero countable ordinal $\alpha$ the space $X_{\mathcal{S}_\alpha}$ (real or complex) has a diagonal standard group of isometries. This result is a consequence of the following general result for surjective isometries between combinatorial spaces $X_{\mathcal{F}}$ and $X_\mathcal{G}$.
	\begin{theorem}
		Let $\mathcal{F}$ and $\mathcal{G}$ be regular families of finite subsets of $\mathbb{N}$. Then the following are equivalent:
		\begin{enumerate}
			\item $T:X_\mathcal{F}\to X_\mathcal{G}$ is a surjective isometry.
			\item $T:X^*_\mathcal{F}\to X^*_\mathcal{G}$ is a surjective isometry.
			\item $Te_i=\varepsilon_i e_{\pi(i)}$ for $\varepsilon_i \in \{\pm1\}$ and some permutation $\pi$ of $\mathbb{N}$ satisfying $\mathcal{G}=\{\pi(F): F\in \mathcal{F}\}$. 
		\end{enumerate}  
	\end{theorem}
	The proof of the above theorem uses the fact, first stated by Gowers \cite{Go-blog}, and proved in \cite{ABC-Illinois}, that 
	$$\ext(B(X_{\mathcal{F}}))= \{ \sum_{i \in F} \varepsilon_i e^*_i : F \in \mathcal{F}^{\text{max}} \text{ and $(\varepsilon_i)_{i\in F}$ is a sequence of signs}\}.$$
	Most recently, Brech and Pina \cite[Theorem 7]{BrechPina-comb} proved a very general result for combinatoral spaces for families on any ordinal $\kappa$. They also prove the remarkable fact for regular families $\mathcal{F}$ and $\mathcal{G}$ of $\mathbb{N}$, the Banach space $X_\mathcal{F}$ and $X_\mathcal{G}$ are isometric if and only if $X_\mathcal{F}^*$ and $X_\mathcal{G}^*$ are isometric  if and only if $\mathcal{F}=\mathcal{G}$. In \cite{BFT-PAMS} they pose the following.
	
	\begin{problem}
		For which combinatorial Banach spaces can we explicitly describe the group of surjective isometries?
	\end{problem}
	
	The following linear extension problem restricted to combinatorial Banach spaces is open.
	
	\begin{problem}
		Let $\mathcal{F}$ be a regular family and $X$ and Banach space. Suppose that $V:S_{X_\mathcal{F}}\to S_X$ is a surjective (not necessarily linear) isometry. Does there exist a linear surjective isometry $U:X_\mathcal{F}\to X$ so the $U|_{S_\mathcal{F}}=V$?
	\end{problem}
	
	This problem is open even in the case $\mathcal{F}=\mathcal{S}_1$. Let us make a few remarks on this problem. A space having the desired property is said to have the Mazur-Ulam property. Famously, Mazur and Ulam 
	proved that a surjective isometry between Banach spaces that fixes $0$ is a linear isometry. There is a vast literature verifying that various spaces have the Mazur-Ulam property \cite{fleming-isometries}. In particular, all classical sequence spaces have the property \cite{Ding-Tingc0, FangWang-ellp, KadetMartin}. As for non-classical sequence spaces, Tan showed that surjective isometries on spheres of Tsirelson's space \cite{Tan-Tsirelson} and James' space (in one of the natural norms) \cite{Tan-James} can be extended to linear isometries on the whole space. Indeed, she showed that each surjective isometry on the sphere has the same form as a surjective isometry on the space. Theorems that do not assume the linearity of the isometry often require a distinct analysis from the linear case. All of these theorems contribute to the well-known linear extension problem (or Tingley's problem \cite{Tingley}) which asks: 
	
	\begin{problem}
		If $X, Y$ are Banach spaces and $V:S_X \to S_Y$ is a surjective isometry is there a linear surjective isometry $U: X\to Y$ so that $U|_{S_X}=V$?
	\end{problem}
	
	This is a major open problem that has only recently been settled in dimension 2 \cite{Banakh-Tingley}. In the final section of this current paper, we prove the Mazur-Ulam Theorem and show that $c_0$ has the Mazur-Ulam property.

	\section{Isometries on Tsirelson-type spaces}
	
	In this section, we will characterize the surjective isometries on Tsirelson's space and some of its variants. Tsirelson's space \cite{Ts-space} is the first example of a space containing no isomorphic copy of $\ell_p$ for $1\leqslant p<\infty$ or $c_0$. The space now called Tsirelson's space is the dual of the original space and was first presented in the form which we present by Figiel and Johnson \cite{FJ-Tsirelson}. Tsirelson's space $T$ is a sequence space defined as the completion of $c_{00}$ with respect to an inductively defined norm and, remarkably, the norm satisfies an implicit formula. The space $T$ is reflexive and $(e_i)$ is a $1$-unconditional basis. Although $T$ has been extensively studied in Banach space theory and the method of construction has led to numerous breakthroughs, it is fair to say that the construction is not well understood outside of Banach space theory. Nevertheless, Tsirelson's example is now almost 50 years old and a seminal work of the isomorphic theory of Banach spaces. Fortunately, we do not need many of the isomorphic properties of $T$ to study $\Isom(T)$. We believe that the reader unfamiliar with this construction will still be able to follow the basic argument.
	
	It is interesting to note that $T$ and the variants we consider do not have a standard or diagonal standard group of isometries, but instead, the complexity of the isometry group lies between these two classes. 
	
	For $x\in c_{00}$ and $E\subset \mathbb{N}$ let $Ex$ be the projection of the vector $x$ onto the coordinates of $E$.
	
	We define a Banach sequence space denoted $T[\theta,\mathcal{F}]$ for a scalar $\theta\in (0,1)$ and a regular family $\mathcal{F}$ and $T:=T[\frac{1}{2},\mathcal{S}_1]$. We chose this presentation because it allows us to ask more general questions which extend the best known results. In addition, the construction of the space $T[\theta,\mathcal{F}]$ is not more complicated than that of $T$.
	
	Let $\|\cdot\|_0$ denote the supremum norm on $c_{00}$. Assuming that $\|\cdot\|_n$ has been defined let
	$$\|x\|_{n+1}= \max\{ \|x\|_n,  \sup \{ \theta \sum_{i=1}^d \|E_ix\|_n :E_1 <\cdots < E_d \text{ in }\mathbb{N}, (\min E_{i})_{i=1}^d \in \mathcal{F}\}\}$$
	Let $\|x\|_{\theta,\mathcal{F}}=\sup_{n \in \mathbb{N}}\|x\|_n$ and $T[\theta,\mathcal{F}]$ be the Banach space equipped with this norm that is the completion of $c_{00}$. A proof by induction shows that the norm satisfies the following implicit formula for $x \in T[\theta,\mathcal{F}]$: 
	$$\|x\|_{\theta,\mathcal{F}}= \max\{|x\|_{\infty} ,\sup\{\theta \sum_{i=1}^d \|E_ix\|_{\theta,\mathcal{F}} : E_1 <\cdots < E_d \text{ in }\mathbb{N}, \,(\min E_i)_{i=1}^d \in \mathcal{F}\}\}$$  
	It follows easily from the definition that if $(x_i)_{i=1}^d$ is a block sequence in $T[\theta,\mathcal{F}]$ with $(\min\text{supp}\,x_i)_{i=1}^d\in \mathcal{F}$ we have  
	\begin{equation}
	\label{eq55}
	\|\sum_{i=1}^n x_i\|_{\theta,\mathcal{F}} \geqslant \theta\sum_{i=1}^n\|x_i\|_{\theta,\mathcal{F}}. 
	\end{equation}
	Restricted to Tsirelson's space $T[\frac{1}{2},\mathcal{S}_1]$ the above states that if $d\leqslant x_1<\cdots <x_d$ then 
	\begin{equation}
	\|\sum_{i=1}^d x_i\|_{\frac{1}{2},\mathcal{S}_1} \geqslant \frac{1}{2}\sum_{i=1}^d\|x_i\|_{\frac{1}{2},\mathcal{S}_1}   \label{lower}
	\end{equation}
	For more information on the spaces $T[\theta,\mathcal{F}]$ see \cite{ADel-TAMS,ATo-book,ATol-Memoirs}.
	
	Our main theorem concerns the spaces $T[\frac{1}{n},\mathcal{S}_1]$ for $n \geqslant 2$. The result for $n=2$ is due to Beauzamy and Casazza and appears in \cite[Theorem III.8]{CS-book}.

	\begin{theorem}
		Let $n \in \nn$ with $n\geqslant 2$. Then $U \in \Isom(T[\frac{1}{n},\mathcal{S}_1])$ if and only if there is a permutation $\pi$ on $\{1,\cdots, n\}$ and $(\varepsilon_i)_{i=1}^\infty\subset \{\pm1\}$ so that $Ue_i=\varepsilon_i e_{\pi(i)}$ for $1\leqslant i\leqslant n$ and $Ue_i = \varepsilon_i e_i$ for $i \geqslant n+1$.
		\label{Tisom}
	\end{theorem}

	\begin{proof}
		Let $U \in \Isom(T[\frac{1}{n},\mathcal{S}_1])$. As $c_{00}$ is dense in $T[\frac{1}{n},\mathcal{S}_1]$, we will make the convention that all vectors considered in the proof are in $c_{00}$. This will allow us to suppress many $\varepsilon$ approximations throughout the proof. For each $j \in \mathbb{N}$ let $Ue_j= \sum_{i=1}^\infty a^j_i e_i $.
		
		Fix $j \in \{1,\dots,n\}$ and let $x= \sum_{i=n+1}^\infty a^j_i e_i $. Since the sequence $(Ue_j)_{j=1}^\infty$ is coordinate-wise null we can find indices $j_1<\cdots <j_n$ with $x< Ue_{j_1} < \cdots < Ue_{j_n}$ (recall that by our convention $Ue_j \in c_{00}$). Therefore
		\begin{equation}
		\begin{split}
		1  & = \|e_j + \sum_{i=1}^n e_{j_i}\| = \|\sum_{i=1}^n a^j_i e_i + x + \sum_{i=1}^n U e_{j_i}\|\geqslant  \frac{1}{n} (\|x\| + n) 
		\end{split}
		\end{equation}
		The first equality follows from computing the norm and final inequality follows from \eqref{eq55}. Thus $\|x\|=0$. Consequently $U[e_1, \ldots, e_n]=[e_1, \ldots, e_n]$. We will show that for $j \in \{1,\ldots ,n\}$, $Ue_j=\pm e_{\pi(j)}$ for some permutation $\pi$ of $\{1,\ldots , n\}$. By definition of the norm, for each such $j$ there at is least one index $\pi(j)$ so that $|a^j_{\pi(j)}|=1$. Also by definition if $i\not=j$ in $\{1,\ldots, n\}$ we have $1=\|e_j\pm e_i\|\geqslant |a_\ell^i\pm a_\ell^j|$ for each $1 \leqslant \ell \leqslant n $. Therefore $|a^j_\ell|=0$ if $\ell\not=\pi(j)$. Thus $\pi$ is the desired permutation. 
		
		We now prove that $Ue_{n+1}=\pm e_{n+1}$. This will be the base case of an induction, however, the proof of the inductive step is very similar. Let $x=\sum_{i=n+3}^\infty a_i^{n+1}e_i$. As before we claim that $\|x\|=0$. Find $j_1<\cdots<j_{n+2}$ so that $x< Ue_{j_1} < \cdots < Ue_{j_{n+2}}$. Then 
		$$\frac{n+2}{n}= \|e_{n+1}+ \sum_{i=1}^{n+2} e_{j_i}\|\geqslant \frac{1}{n}(\|x\| + n+2)$$
		Thus $\|x\|=0$. In the above we used \eqref{eq55} since $n+3 \leqslant \supp \, x$. For $1 \leqslant j \leqslant n$ we have 
		$$1=\|e_j\pm e_{n+1}\|= \|e_{\pi(j)}\pm Ue_{n+1}\| \leqslant |1\pm a^{n+1}_{\pi(j)}|$$
		Therefore $a^{n+1}_{\pi(j)}=0$ for all such $j$. Whence we have that $Ue_{n+1}=a^{n+1}_{n+1} e_{n+1}+a^{n+1}_{n+2} e_{n+2}$. We claim that $a^{n+1}_{n+2}=0$. Find $n+2<j_1<\cdots<j_{n+1}$ so that $n+2< Ue_{j_1} < \cdots < Ue_{j_{n+1}}$. Then 
		$$\frac{n+1}{n}=\|e_{n+1} + \sum_{i=1}^{n+1}e_{j_i}\| = \|a^{n+1}_{n+1} e_{n+1}+a^{n+1}_{n+2}e_{n+2} + \sum_{i=1}^{n+1}U e_{j_i}\|\geqslant \frac{1}{n}(|a^{n+1}_{n+2}| + n+1)$$
		Thus $|a^{n+1}_{n+2}|=0$. It follows that $Ue_{n+1}=\pm e_{n+1}$. 
		
		We now proceed to the inductive step. Assume that for some $m >n+1$ and for all $n<k\leqslant n+1$, $Ue_k=\pm e_k$. Then by slightly modifying the previous arguments we can see that $Ue_m= a^m_me_m + a^m_{m+1}e_{m+1} + x$ where $x= \sum_{k=m+2}^\infty a^m_{k}e_k$. Find $j_1<\cdots <j_{m+1}$ with $x<Ue_{j_1}< \cdots < Ue_{j_{m+1}}$. Then again we have $\|x\|=0$ since
		$$\frac{m+1}{n}=\|e_m +\sum_{i=1}^{m+1} e_{j_i}\|\geqslant \frac{1}{n}(\|x\|+m+1).$$
		The same argument as in the base case shows $Ue_m=\pm e_m$. This finishes the proof that the surjective isometry has the desired form. 
		
		It remains to show that for all such maps are indeed isometries. Let $\pi$ be permutation of $\{1,\cdots , n\}$. It suffices to show that for each $x \in c_{00}$ with $n< \supp \, x$ we have
		$$\|\sum_{i=1}^n a_i e_i +x \|= \|\sum_{i=1}^n a_{i} e_{\pi(i)} +x\|$$
		Let $y_1 = \sum_{i=1}^n a_i e_i +x$ and $y_2 = \sum_{i=1}^n a_{i} e_{\pi(i)} +x$. We will show that $\|y_1\|\leqslant \|y_2\|$. The reverse direction is the same. The vectors $y_1$ and $y_2$ have the same supremum norm. So we assume that norm of $y_1$ given by sets $2\leqslant d \leqslant E_1<\cdots < E_d$ in the sense that 
		\begin{equation}
		\|y_1\|= \frac{1}{n}\sum_{i=1}^d\|E_i y_1\|.
		\label{breaking}
		\end{equation}
		If $n<E_1$ then $\|E_iy_1\|=\|E_iy_2\|$ for $i$ and so $\|y_1\|\leqslant \|y_2\|$. If $\min E_1 <n$. Then (\ref{breaking}) cannot hold since $\frac{1}{n}\sum_{i=1}^d\|E_iy_1\|\leqslant \frac{d}{n}\|y_1\|$ and $d\leqslant \min E_1<n$. If $\min E_1=n$ then $d=n$ and from (\ref{breaking}) $\|y_1\|=\|E_iy_1\|$ for each $1\leqslant i \leqslant n$. However for $i\geqslant 2$ we have $\|E_iy_1\|=\|E_iy_2\|$. Thus for each $2\leqslant i \leqslant n$
		$$\|E_1y_1\|=\|y_1\|=\| E_i y_1\|=\|E_iy_2\|\leqslant \|y_2\|.$$
		Consequently, $\|y_1\|\leqslant \|y_2\|$ as desired. 
	\end{proof}

	In the proof of the above theorem, the assumption that $U$ was surjective was not used. Therefore we have the following corollary.
	
	\begin{corollary}
		Let $n \in \nn$ with $n \geqslant 2$. Then every linear isometry on $T[\frac{1}{n},\mathcal{S}_1]$ is surjective. 
	\end{corollary}
	
	We do not have a description of the surjective isometries on $T[\theta, \mathcal{F}]$ for many other choices of $\theta \in (0,1)$ and regular families $\mathcal{F}$. The case of $\theta =\frac{1}{2}$ and $\mathcal{F}=\mathcal{S}_{\alpha}$ for some countable ordinal $\alpha>1$ is especially interesting. 
	
	A further class of spaces that generalizes Tsirelson's space are the mixed-Tsirelson spaces denoted $T[(\theta_n,\mathcal{F}_n)_{n=1}^\infty]$ for decreasing sequence $(\theta_n)$ and regular families $(\mathcal{F}_n)$. These space are extensive studied in \cite{ADel-TAMS,ATo-book,ATol-Memoirs}. The norm of $T[(\theta_n,\mathcal{F}_n)_{n=1}^\infty]$ satisfies the following implicit equation. 
	$$\|x\|= \max\{\|x\|_{\infty},\sup\{\theta_n \sum_{i=1}^d \|E_ix\| : E_1 <\cdots < E_d \text{ intervals in }\mathbb{N}, \,(\min E_i)_{i=1}^d \in \mathcal{F}_n\}\}.$$
	Tsirelson's space $T[\frac{1}{2},\mathcal{S}_1]$ is naturally isometric to the mixed Tsirelson space $T[(\frac{1}{2^n},\mathcal{S}_n)_{n=1}^\infty]$. 
	Let $\mathcal{A}_n$ be the finite subsets of $\mathbb{N}$ with cardinality less than $n$ and $\theta_n=1/\log_2(n+1)$. Then $T[(\theta_n,\mathcal{A}_n)_{n=1}^\infty]$ is the famous space of Th. Schlumprecht \cite{Sc-Israel}; the first known arbitarily distortable Banach space and the precursor to the hereditarily indecomposable space of Gowers and Maurey \cite{GM}. The geometry of the mixed Tsirelson spaces is, of course, very sensitive to the parameters used to define the space. For example for $\theta_1\not=\theta_2$ in $(0,1)$ the space $T[\theta_1,\mathcal{S}_1]$ is not isomorphic to a closed subspace of $T[\theta_1,\mathcal{S}_1]$. As far we know, Theorem \ref{Tisom} is the only characterization of the isometries on Tsirelson-type spaces or mixed Tsirelon space.  
	
	\begin{problem}
		Characterize the isometries on a class of mixed Tsirelson spaces.  
	\end{problem}
	
	For particularly interesting classes of mixed Tsirelson space see \cite{AO-Israel} or \cite{MP-BLMS}.

	\section{Spaces with Small Groups of Isometries}
	
	In this section, we present some spaces which have a trivial group of surjective isometries and highlight some results concerning other groups of isometries that are group isomorphic to the group of surjective isometries. 
	
	The first natural question is whether a space with only trivial isometries exists. Our starting point is the following example due to A. Pe\l czynski (see Davis' paper \cite{Davis-Isom}). 
	
	\subsection{A $C(K)$ space with a trivial group of isometries} Using the Banach-Stone theorem, it will suffice to find an infinite compact metric space so that the only surjective automorphism in the identity. A Hausdorff space so that the only surjective automorphism is the identity is called rigid and the first such space was constructed by Kuratowski in 1926. 
	To construct a metric example one may take $K$ to be a compact, locally connected metric space that does not contain a simple closed curve (i.e. dendrite) so that for each $n\geqslant 3$, $K$ contains is a unique cut point of degree\footnote{The degree of the cut-point is the number of disjoint components the space is broken into upon deletion.} $n$ and such that the set of all such cut points is dense in $K$.

	Constructing such a space goes as follows. Start with a closed interval $I_1$ of length $1$ in the plane. Place a closed line segment $I_2$ of length $1/2$ so that one end is at the midpoint of the $I_1$, call this midpoint $x_3$. Clearly $x_3$ is a degree $3$ cut point. Now pick one of the three line segments that remain when removing $x_3$ and call it $I_3$. Place two different line segments of length $1/3$ that have endpoints at the midpoint $x_4$ of $I_3$. Then $x_4$ is a cut point of degree $4$. Call this the second iteration. Choose $I_5$ to be one of the other two intervals created when removing $x_3$, and repeat this procedure: At the $n^{th}$ iteration place $n$ line segments of length $1/(n+1)$ having endpoints at the midpoint $x_{n+2}$ of the interval $I_{n+1}$. The resulting space $K$ has, for each $n>2$, a unique cut-point of degree $n$, and cut points are dense in $K$ as desired. \\

	In 1971, Davis \cite{Davis-Isom} showed that the following equivalent norm on $\ell_2$ has a trivial group of isometries. Moreover, there are no non-surjective isometries in this renorming.
	
	\subsection{Davis' renorming of $\ell_2$}
	Let $X_2=\mathbb{R}\oplus \ell_2$ be a vector space, let $e_0$ be a unit vector in $\mathbb{R}$, and $(e_i)_{i=1}^\infty$ be the unit vector basis of $\ell_2$. For each $n \in \mathbb{N}$ define the set $\Lambda_n=\{(0,t e_n) : -\frac{1}{2n+1}\leqslant t \leqslant \frac{1}{2n}\}$ in $\{0\} \oplus \ell_2 $ and let $F = \overline{\mbox{co}}\{\Lambda_n:n\in \mathbb{N}\}$, where $\overline{\mbox{co}}(X)$ denotes the closed convex hull of $X$. We define the norm on $X_2$ by defining unit ball of $X_2$. Set
	$$B_{X_2}=\overline{\mbox{co}}\{F+(0,e_0),-F-(0,e_0),B_{\ell_2\oplus \{0\}}\}.$$
	
	Davis' method is applicable to the other classical spaces and the proof of the above theorem is quite geometric in nature. This result naturally led to the question of whether every Banach space can be equivalently renormed so that the group of isomorphisms is trivial. 
	
	In 1986, S. Bellenot \cite{Bellenot-trivial} showed that every separable Banach space has an equivalent norm with only trivial surjective isometries. A year later,  K. Jarosz \cite{J-Isom} proved that every Banach space has such an equivalent renorming. These results imply that it is not possible to isomorphically distinguish Banach spaces by their isometry groups. The constructions of Bellenot and Jarosz are quite involved and so we do not present them here. In the case of $c_0$, however, we can present a renorming having a trivial group of isometries that is elegant and short. The proof is due to Ferenczi and Rosendal \cite{FR-Extracta} and is a special case of a much more general theorem. 
	
	\subsection{A renorming of $c_0$}
	Consider a graph $\Gamma$ on $\mathbb{N}$ having a trivial automorphism group. Using a construction similar to that of the dendrite in the previous section, it is easy to construct a graph so that each vertex labeled $n$ has valence number is exactly $n$. Let $d_\Gamma(m,n)$ be the graph metric and note that for $\pi \in \text{SYM}(\mathbb{N})$ we have
	\begin{equation}
	\pi=\Id \iff d_\Gamma(m,n)=d_\Gamma(\pi(m),\pi(n))\text{ for all }m,n \in \mathbb{N}.
	\label{id}
	\end{equation}
	
	Define the following equivalent norm on $c_0$
	\begin{equation}
	\begin{split}
	\|(a_i)\|_{\Gamma}& = \max_{n,m \in \nn, m\not=n} 
	\{|a_n+\frac{a_m}{1+2d_{\Gamma(m,n)}}|,|a_n-\frac{a_m}{2+2d_{\Gamma(m,n)}}|, |a_n|\}\\
	\end{split}
	\end{equation}

	\begin{lemma}
		$\ext(B_{c_0,\|\cdot\|_{\Gamma}})=\{\pm e_n: n\in \nn\}$.\label{extreme}
	\end{lemma}     
	
	\begin{proof}
		First note that the definition of the norm easily yields that $x=e_n$ for some $n\in \mathbb{N}$ if and only if $\|x\|_\Gamma=\|x\|_\infty$.
		
		Let $n\in \mathbb{N}$ and suppose that $\lambda\in(0,1)$ and 
		$e_n=\lambda(\sum_i c_ie_i)+(1-\lambda)(\sum_i b_i e_i)$ with $\|\sum_i c_ie_i\|,\|\sum_i b_ie_i\|\leqslant 1$. Since $|c_i|,|b_i|\leqslant 1$ and $\lambda c_n+(1-\lambda)b_n=1$ we have $c_n=b_n=1$. From the definition of the norm if follows that $c_i=b_i=0$ for $i\not=n$. Thus $e_n\in \ext(B_{c_0,\|\cdot\|_{\Gamma}})$.
		
		In the other case, let $\sum_i a_i e_i\in c_0$ and assume that $\|\sum_i a_i e_i\|_\infty <1$. Find $\varepsilon>0$ with $\|\sum_i a_i e_i\|_\infty<1-\varepsilon$ and $m \in \mathbb{N}$ so that $|a_m|<\varepsilon/2$. Then clearly
		$$\sum_i a_i e_i = \frac{1}{2}(\underbrace{\sum_i a_i e_i + \frac{\varepsilon}{2}e_m}_{=x})+\frac{1}{2}(\underbrace{\sum_i a_i e_i - \frac{\varepsilon}{2}e_m}_{= y})$$
		Calculating with the norm shows that $\|x\|,\|y\|\leqslant 1$. This is the desired result.
	\end{proof}
	
	\begin{theorem}
		The Banach space $(c_0,\|\cdot\|_{\Gamma})$ has a trivial group of isometries. \label{c0}
	\end{theorem}
	
	\begin{proof}
		Let $T\in \Isom(c_0,\|\cdot\|_{\Gamma})$. By Lemma \ref{extreme} for some $\pi\in \text{SYM}(\mathbb{N})$ we have $Te_n=\pm e_{\pi(n)}$. For $m\not= n$ we have 
		$$1+\frac{1}{1+2d_\Gamma(m,n)}\|e_n+e_m\|\not=\|e_n-e_m\|=1+\frac{1}{2+2d_\Gamma(m,n)}$$
		In fact, for $i\not=j$, the parity of the denomintor yields that
		$$\|e_n+e_m\|\not=\|e_j-e_i\|$$
		Therefore either $Te_n=e_{\pi(n)}$ or $Te_n=-e_{\pi(n)}$. Assume the former. Then 
		$$\|e_n+e_m\|=\|e_{\pi(n)}+e_{\pi(m)}\|$$
		This implies that $d_\Gamma(m,n)=d_\Gamma(\pi(m),\pi(n))$. By (\ref{id}) we conclude that $\pi=\Id$, as desired.
	\end{proof}
	
	\subsection{Displaying groups on a Banach space $X$} Let $GL(X)$ be the group of surjective linear isomorphisms of $X$. A subgroup $G$ of $GL(X)$ is bounded if $\sup_{T\in G}\|T\|<\infty$. Given a bounded subgroup $G$ of $GL(X)$ the following is an equivalent norm on $X$ so that each element of $G$ is an isometry:
	$$\|x\|_G=\sup \{\|Tx\|:T\in G\}$$
	That is, $G\leqslant \Isom(X,\|\cdot\|_G)$.
	
	For any group $G$ (not necessarily a subgroup of $GL(X)$), it is natural to ask if there an equivalent norm $\|\cdot\|$ on $X$ so that $\Isom(X,\|\cdot\|)$ group isomorphic to $G$. If this is the case, we say $G$ is \textit{displayable} on the Banach space $X$. We will distinguish the term displayable and representable in the following way. Classically a representation of a group $G$ on a space $X$ is a homomorphism of $G$ into $GL(X)$. A group has a faithful representation if the representation is injective. In this case, $G$ is group isomorphic to
	a subgroup of $GL(X)$.
	If the group $G$ has a faithful representation by a collection of uniformly bounded isometries on a Banach space then $G$ is said to be \textit{representable}\footnote{The term displayable was introduced by V. Ferenczi and C. Rosendal in \cite{FR-Extracta}. Before this paper, the term representable was used in place of displayable. As we explained, however, since this usage does not cohere with the classical notion of an injective representation, it's prudent to adopt the term displayable.}.  
	A trivial requirement for a group $G$ to be representable on $X$ is for it to have a normal subgroup of order $2$ isomorphic (in order to have a subgroup isomorphic to $\{\pm \Id\}$); that is, $G$ contains a non-trivial central involution. J. Stern \cite{Stern} proved that for any group $G$ containing a non-trivial central involution, there is a real Hilbert space $H$ such that $G$ is displayable on $H$.    
	
	One may also consider the topology when displaying groups on a Banach space $X$. Indeed, since the group operations on $\Isom(X)$ are continuous in the topology of pointwise convergence (i.e. the strong operator topology), $\Isom(X)$ is a topological group that is a Polish (i.e. a complete, separably metrizable topological space) group when $X$ is separable. Therefore we say that a group is topologically displayable if it is topologically isomorphic to the isometry group $\Isom(X,\|\cdot\|)$ for some equivalent norm $\|\cdot\|$ on $X$.
	
	As noted, Theorem \ref{c0} is a particular case of the following much more general theorem of Ferenczi and Rosendal \cite[Theorem 22]{FR-Extracta}. 
	
	\begin{theorem}
		Any closed subgroup $G$ of $\text{SYM}(\nn)$ with a non-trivial central involution is  topologically displayable on $c_0$, $C([0,1])$, $\ell_p$ and $L_p$, for $1 < p < \infty$.
	\end{theorem}
	
	The proof of the above theorem in the case of $c_0$ is similar to the proof Theorem \ref{c0}. Indeed, one only has to consider graphs $\Gamma$ so that the automorphism group of $\Gamma$ is topologically isomorphic to the given closed subgroup $G$ of $\text{SYM}(\mathbb{N})$. In the case of  $\ell_p$, the authors generalize the proof of Bellenot showing that every separable Banach space can be renormed to have a trivial group of isometries. 
	
	We warmly recommend the papers of Ferenczi and Rosendal \cite{FR-Duke, FR-Extracta}, and Ferenczi and Galego \cite{FGalego-TAMS} in which the authors prove many striking theorems involving displays of groups on various Banach spaces. In particular, Ferenczi and Galego prove the following extension of Bellenot's result \cite[Theorem 11]{FGalego-TAMS}.
	
	\begin{theorem}
		Let $X$ be a real separable Banach space and $G$ be a finite subgroup of $GL(X)$ so that $\{-\Id, \Id\}\subset G$. Then there is an equivalent norm $\|\cdot\|$ on $X$ so that $\Isom(X,\|\cdot\|)=G$.\label{FG11}
	\end{theorem}
	
	Trivially the finite group $G$ in the above theorem is representable on $X$, the conclusion of the theorem is that it is also displayable. If we do not assume that the finite group $G$ is a group of isometries, it is natural to ask whether $\{-1,1\}\times G$ can be displayed on $X$. The following result \cite[Theorem 17]{FGalego-TAMS} gives a positive answer to this question. We include the short proof which follows from Theorem \ref{FG11}. 
	
	\begin{theorem}
		Suppose that $G$ is a finite group and $X$ is a real separable Banach space with $\dim X \geqslant |G|$. Then $\{-1,1\}\times G$ is displayable on $X$. \label{finitegps}
	\end{theorem}
	
	\begin{proof}
		Fix $G$ and $X$ as in the hypothesis and consider the finite dimensional Banach space $\ell_2(G)$ with unit vector basis $(e_g)_{g\in G}$. The map
		$$\{- 1,1\}\times G \ni (\varepsilon,g) \mapsto T_{\varepsilon,g}\in \Isom(\ell_2(G))   \text{ where }T_{\varepsilon,g} (\sum_{h \in G}a_h e_h)=\varepsilon (\sum_{h \in G}a_h e_{hg})$$
		is a group isomorphism. Moreover since $\dim X \geqslant |G|$, $X$ is isomorphic to $\ell_2(G)\oplus_2 Y$ for some $Y$ (note that we may take $Y=X$ when $X$ is infinite dimensional). This isomorphism yields a renorming of $X$. For $(\varepsilon,g)\in G$, let $A_{\varepsilon,g}\in \Isom(\ell_2(G)\oplus_2 Y)$ be defined by $A_{\varepsilon,g}(t,x)=(T_{\varepsilon,g}(t),\varepsilon)$. Then the collection of all such operators are a finite group of isometries containing $\{\pm \Id\}$. The conclusion therefore follows from Theorem \ref{FG11}. 
	\end{proof}

	The papers \cite{FGalego-TAMS, FR-Extracta} contain some extensions of the above result for certain countable and uncountable groups and for various Banach spaces. For example, \cite[Theorem 31]{FR-Extracta}:
	
	\begin{theorem}
		Let $G$  be a Polish group. Then there is a separable Banach space $X$ so that $\{-1,1\}\times G$ is displayable on $X$.
	\end{theorem} 
	
	It is natural to ask whether for any group $G$ and for any real Banach space $X$ if $\{-1,1\}\times G$ is displayable on $X$. This turns out to be false for countable groups. The known counterexamples come from a class of spaces with very rigid structure called hereditarily indecomposable spaces. A space $X$ is called hereditarily indecomposable (HI) if no infinite-dimensional closed subspace $Y$ is the topological direct sum of two of its infinite-dimensional closed subspaces. The first example of such a space is due to Gowers and Maurey \cite{GM}, a paper that is now a seminal work in the isomorphic theory of infinite-dimensional Banach spaces. It is an easy fact that such spaces cannot contain unconditional basic sequences. The following holds \cite[Proposition 7.8]{FR-Duke}:
	
	\begin{theorem}
		There is a real separable HI space $X$ (without a Schauder basis) so that for any norm on $X$, the group of isometries on $X$ is either finite or of cardinality $2^{\aleph_0}$. \label{split}
	\end{theorem}
	
	Notice that is it not possible to exclude the second possibility. Indeed, for any $n \in \mathbb{N}$ and Banach space $X$, the space $\ell_2^n\otimes X$ is isomorphic to $X$ and therefore an equivalent renorming of $X$. The isometry group on $\ell_2^n\otimes X$ will at least contain a copy of the unitary group for dimension $n$, which has cardinality $2^{\aleph_0}$. 
	
	Ferenczi and Rosendal deduce Theorem \ref{split} and many other results after first showing that each isometry on a Banach space having \textit{few operators}\footnote{A space $X$ has few operators if every bounded linear operator $T:X \to X$ has he form $T=\lambda I+S$ for $\lambda\in \mathbb{R}$ and $S$ a strictly singular operator. Complex hereditarily indecomposable Banach spaces \cite{GM} have this property as well as many real HI spaces including the original space of Gowers and Maurey.} and not containing an infinite-dimensional subspace with unconditional basis has a very rigid form. The following is a consequence of \cite[Theorem 6.6]{FR-Duke}:
	
	\begin{theorem}
		Let $X$ be a Banach space not containing an unconditional basic sequence and having few operators. Then $U\in\Isom (X)$ is of the form $U=\varepsilon \Id+F$ where $|\varepsilon|=1$ and $F$ is a finite rank operator; that is, a finite rank perturbation of the identity. 
	\end{theorem}
	
	A subgroup $G$ of $GL(X)$ is SOT-discrete if it is discrete in the topology of pointwise convergence. Finite groups are trivially SOT-discrete, the following question is stated in \cite[Problem 8.4]{AFGR-light}.
	
	\begin{question}
		Does there exist a reflexive space $X$ with an SOT-discrete infinite bounded group $G \leqslant GL(X)$ such that all elements of $G$ are finite rank perturbations of the identity?
	\end{question}
	
	The space $c$ (see \cite[Lemma 4.3]{AFGR-light}) is a non-reflexive example of such as space: Consider the group of isometries of the form $\Id+F$ so that $F$ is a sign change of finitely many coordinates of
	the sequence. This group is discrete because for any such non-trivial operator $U$ of this form,  $\|(U-\Id)(e)\|_\infty=2$, where $e$ is the constant sequence 1.
	We highlight one more problem related to our current presentation. 
	
	\begin{problem}
		Let $X$ be a combinatorial Banach space. Which closed subgroups $G$ of $\text{SYM}(\mathbb{N})$ can be displayed on $X$? 
	\end{problem}
	
	Other than Theorem \ref{split}, to the authors' knowledge, there is no other class of spaces (e.g. combinatorial, Tsirelson-type) for which there is a known restriction on the cardinality of the groups of isometries that are representable (or displayable) on $X$. That is, there is a lot of work to be done in this area.

	While there are many deep and interesting problems concerning which groups can be represented or displayed on various classes of Banach spaces it seems, however, that there is also much more work to be done in characterizing $\Isom(X)$ for different Banach spaces $X$ equipped their standard norms. Indeed such characterizations could be considered the natural first step in studying the groups that are displayable or representable on $X$. 
	
	In the final subsections of this section, we present some spaces that have a trivial group of isometries in their natural norm and pose some related questions. 
	
	First note that if one wishes to find an example of a Banach space whose natural norm has a trivial group of isometries, the space should at minimum not have a $1$-unconditional Schauder basis. The James space $J_2$ is thus a natural example to consider and, historically, the literature starts with this space. 
	
	\subsection{Isometries on James space} The James space $J_2$ of all null sequences of reals with bounded square variation, constructed by R.C. James \cite{Ja-PNAS} in 1951. James' space has played an important role in the development of Banach space theory. Amongst its other interesting properties, $J_2$ has the property that it cannot be isomorphically embedded into a space with an unconditional basis. The space $J_2$ was the first example of a non-reflexive Banach space that is isomorphic to its bi-dual. Moreover, $J_2$ is quasi-reflexive of order $1$, that is $\dim(J_2^{**}/J_2)=1$. A good reference for James space is the book of Albiac and Kalton \cite{AlbK-book}. Several equivalent norms have been considered on $J_2$. The following two norms are the standard norms on $J_2$:
	\begin{enumerate}
		\item Let $x=(a_i)\in c_0$ 
		$$\|x\|_{1}= \frac{1}{\sqrt{2}} \sup \{(\sum_{i=1}^n |a_{p_i}-a_{p_{i+1}}|^2)^{1/2}: 0\leqslant p_1<p_2< \cdots <p_{n+1}, n \in \mathbb{N}\} $$ 
		We assume that $a_0=0$ in the case that $p_1=0$.  
		\item Let $x=(a_i)\in c_0$
		$$\|x\|_2= \frac{1}{\sqrt{2}} \sup \{(\sum_{i=1}^n |a_{p_i}-a_{p_{i+1}}|^2 + |a_{p_{n+1}}-a_{p_1}|)^{1/2}: 0\leqslant p_1<p_2< \cdots <p_{n+1}, n \in \mathbb{N}\} $$
	\end{enumerate}
	$J_2$ is isomorphic to $J_2^{**}$ in both norms but isometric to $J_2^{**}$ in $\|\cdot\|_2$. 
	
	\begin{theorem}
		The space $J$ has a trivial group of isometries in both norms $\|\cdot\|_1$ and $\|\cdot\|_2$.  \label{James}
	\end{theorem}
	
	The case of the norm $\|\cdot\|_1$ is due to S. Bellenot in \cite{Bellenot-James}. He derives this as a consequence of a more general result. Indeed, Bellenot applies his general result to show that by modifying the norm of a certain Orlicz sequence in the same way that James' space is a modification of $\ell_2$, one can construct a space with a non-separable bi-dual, whose natural norm, has a trivial group of isometries. The proof is quite technically demanding and relies on extreme point arguments which do not carry over to the other natural norm. 
	
	A. Sersouri \cite{Sersouri-TAMS} proves the case of Theorem \ref{James} for $J_2$ equipped with the norm $\|\cdot\|_2$ from a general and deep result about quasi-reflexive spaces that are isometric to their biduals. This paper contains many difficult results related to isometries operators on quasi-reflexive Banach spaces.

	As for renormings of $J_2$, Semenov and Skorik \cite{SS-James} were the first to characterize isometries on $J_2$. In addition to providing a non-standard renomorming of $J_2$ for which $J_2$ has trivial isometries, they also show the following renorming of $J_2$:
	
	$$\|(a_i)\|_{3}= \sup \{(\sum_{i=1}^n |a_{p_{2i-1}}-a_{p_{2i}}|^2)^{1/2}: p_1<p_2< \cdots <p_{2n}, n \in \mathbb{N}\} $$
	has the property that $\Isom(J_2,\|\cdot\|_{3})$ is isomorphic to $\mathbb{Z}_2\times \mathbb{Z}_2$ in the real case and $S^1\times \mathbb{Z}_2$ in the complex case.  Note that the operator switching $e_1$ and $e_2$ is an isometry in this norm. Finally, note that since $J_2$ is isomorphic to $J_2 \oplus \ell_2$ it is possible to renorm $J_2$ to have many isometries.

	We conclude this section by listing some spaces and corresponding norms for which we conjecture that the groups of surjective isometries may be trivial.

	\subsection{James-Schreier spaces}
	Let $1 \leqslant p < \infty$, Bird and Laustsen \cite{BirdJamLaus-JOT} (also see \cite{BLZ}) introduced the $p$th James-Schreier spaces $V_p$. These spaces are the completion of $c_{00}$ under the following norm.
	
	$$\|\sum_i a_i e_i\|_{V_p} = \sup \{\bigg(\sum_{i=1}^k |a_{p_i}-a_{p_{i+1}}|^p \bigg)^{1/p}: k \leqslant p_1<\dots < p_{k+1}\}$$
	
	In the case that $p=2$ one has the standard James space with the Schreier-type admissibility condition on the indices $(p_1,\ldots, p_{k+1})$. For each $p$ the $p$th the natural basis is shrinking and every infinite-dimensional subspace contains an isomorphic copy of $c_0$. Therefore, in contrast to $J_2$, $V_p^{**}$ is non-separable. This leads to the following natural question.
	
	\begin{problem}
		Let $1\leqslant p <\infty$. Is there a natural equivalent norm on $V_p$ having only trivial surjective isometries? Which groups are representable or displayable on $V_p$?
	\end{problem}

	\subsection{Jamesifications of Spaces} The method introduced by James to modify $\ell_2$ to construct $J_2$ can be extended to modify any Banach space with a basis. For a Banach space $X$ is a basis basis $(x_n)$ denote by $J(x_n)$ \textit{Jamesification} of $X$. Let $J(x_i)$ be the completion of $c_{00}$ with respect to the following norm. 
	$$\|\sum a_i e_i\|_{J(x_i)}= \sup \{ \|\sum_{i=1}^k(\sum_{j=p(i)}^{q(i)} a_i) x_{p(i)}\|_X: 1 \leqslant p(1)\leqslant q(1)<  \cdots < p(k)\leqslant q(k)\}$$
	This method was defined and studied by Bellenot, Odell, and Haydon \cite{BHO-James}.  In Theorem 4.1 of this work, they proved if $X$ is reflexive then $J(x_i)$ is quasi-reflexive of order-one. This solved a problem of A. Pe\l czy\'nski. 
	If $X=\ell_2$ and $(x_i)$ is the unit vector basis then $J(x_i)$ a norm on James space with the bounded complete basis. Note that the sequence $d_1=e_1$ and $d_i=e_{i}-e_{i-1}$ for $i>1$ forms a Schauder basis for $J(x_i)$ that coincides with the usual (shrinking) basis of $J_2$.
	In \cite{AMS} the authors refer to the above definition as the Jamesification and have a further definition in which they require $q(i)+1=p(i+1)$ call this the conditional Jamesification. These definitions coincide with the basis $(x_i)$ in suppression unconditional.  The following question is open for many reflexive spaces $X$ and unconditional basic sequences $(x_i)$. 
	
	\begin{problem}
		Let $X$ be a reflexive space with an unconditional basis $(x_i)$. Does the space $J(x_i)$ with the given norm have a trivial group of isometries? Which groups of isometries are representable or displayable on $J(x_i)$?
	\end{problem}
	
	\section{Mazur-Ulam Theorem and Linear Extensions}
	
	This final section includes the proof of the Mazur-Ulam Theorem and a theorem of G.G. Ding \cite{Ding-Tingc0} showing that $c_0$ has the Mazur-Ulam property. 
	\begin{theorem}
		Let $X$ and $Y$ real norm spaces. If $T:X \to Y$ is a (not necessarily linear) surjective isometry so that $T(0)=0$. Then $T$ is linear.
	\end{theorem}
	
	\begin{proof}
		Let $T:X\to Y$ be a surjective isometry with $T(0)=0$. 
		The main tools for the argument are following sets: For each $x,y \in X$ and $n \in \nn$ recursively define the nested set $(H_n(x,y))_{n=1}^\infty$. Let $H_1(x,y)=\{z : \|x-z\|=\|y-z\|=\frac{1}{2}\|x-y\|\}$ and having defined $H_n(x,y)$ for some $n \geqslant 1$ let 
		$$H_{n+1}(x,y)=\{z \in H_{n}(x,y):\|z-w\|\leqslant \frac{1}{2}\delta(H_{n}(x,y))\}$$
		where $\delta(H_{n}(x,y))$ is the diameter of the set $H_{n}(x,y)$. Straightforward inductions prove each of the following:
		\begin{enumerate}
			\item For $x,y\in X$ we have $\{\frac{1}{2}(x+y)\}= \cap_{n=1}^\infty H_n(x,y)$.
			\item For $x,y \in X$ and $n \in \nn$ we have $T(H_n(x,y))=H_n(Tx,Ty)$.
		\end{enumerate}
		Combining the (1) and (2) we may conclude the $T(\frac{1}{2}(x+y))= \frac{1}{2}(Tx+Ty)$. That is, $T$ sends centers to centers. From this point on, showing that $T$ is linear is a standard exercise. Let $x,y \in X$. Then 
		$$T(x+y)=\frac{1}{2}T(2x) + \frac{1}{2}T(2y)=\frac{1}{2}(T(2x) +T(0)) + \frac{1}{2}(T(2y)+T(0))= T(x) +T(y).$$
		Therefore $T(nx)=nT(x)$ for $n\in \mathbb{N}$ and since $T(0)=T(x)+T(-x)$ we have $T(-x)=-T(x)$. Thus $T(kx)=kT(x)$ for $k \in \mathbb{Z}$. Whence $T(x)=kT(x/k)$ for $k \in \mathbb{Z}$, $k\not=0$. Thus $T(qx)=qT(x)$ for $q \in \mathbb{Q}$ and linearity follows from the continuity of $T$.\end{proof}
	
	A real Banach space $X$ has the Mazur-Ulam property if for each $Y$ every isometry $V:S_X\to S_Y$ is the restriction of a surjective linear isometry from $X \to Y$. As we previously mentioned, Tingley's problem asking if every real Banach space has the Mazur-Ulam property.  Note that the corresponding question for complex spaces has a trivial counterexample since the conjugate function $f(z)=\overline{z}$ defined on the sphere of $\mathbb{C}$ is an isometry but is not the restriction of a $\mathbb{C}$-linear map on $\mathbb{C}$. 
	
	The following is a special case of a more general result of Ding \cite{Ding-Tingc0}, however, the proof is almost exactly the same.
	We choose to include the proof of Ding's Theorem in this survey since it gives a blueprint towards showing that a given Banach space has the Mazur-Ulam property. We also felt that it was somewhat difficult to extract a self-contained proof of this precise theorem from the existing literature.
	
	\begin{theorem}
		Let $X$ be Banach space and  $V:S_{c_0} \to S_X$ be surjective isometry. Then there is a linear surjective isometry $U:c_0 \to X$ so that $U|_{S_{c_0}}=V$.
	\end{theorem}
	
	We need the following well-known lemma (see \cite[Lemma 2]{Ding-Tingc0}.
	
	\begin{lemma}
		Let $X$ be a normed space and $(x_i)_{i=1}^n\subset S_X$ so that for  $(\varepsilon_i)_{i=1}^n \subset \{\pm 1\}$, $\|\sum_{i=1}^n \varepsilon_i x_i\|=1$. Then for any scalars $(b_i)_{i=1}^n $, $\|\sum_{i=1}^n b_i x_i\|=\sup_{1\leqslant i\leqslant n} |b_i|.$
		\label{allscalars}
	\end{lemma}
	
	\begin{proof}
		We claim that for $x \in S_{c_0}$ we have $V(-x)=-V(x)$. We first prove that $V(-e_n)=-V(e_n)$ for each $n \in \mathbb{N}$. Let $n,m \in \mathbb{N}$ and set
		$$V(x_n)=-V(e_n),V(x_m)=-V(e_m), V(y_n)=-V(-e_n),\mbox{ and }V(y_m)=-V(-e_m).$$
		Let $d$ be the metric associated to the norm. Then
		$$\|x_n-e_n\|=d(V(x_n),V(e_n))=d(-V(e_n),V(e_n))=2.$$
		Likewise $\|x_m-e_m\|=2$. It follows that $x_n(n)=x_m(m)=-1$. Then 
		$$\|x_n-x_m\|=d(-V(e_n),-V(e_m))=d(V(e_n),V(e_m))=\|e_n-e_m\|=1.$$
		Thus $x_n(m)\leqslant 0$. However,
		$$\|y_m+e_m\|=d(y_m,-e_m)=d(-V(-e_m),V(-e_m))=2$$
		and 
		$$\|x_n-y_m\|=d(-V(e_n),-V(-e_m))=\|x_n+x_m\|=1.$$
		Thus $y_m(m)=1$ and $x_n(m)\geqslant 0$. It follows that $x_n=-e_n$.
		
		Let $x\in c_0$ and $y\in c_0$ with $V(y)=-V(x)$. We will show that $y=-x$. We have 
		$$\|e_m+y\|=\|e_m-x\|\mbox{ and }\|e_m-y\|=\|e_m+x\|$$
		Thus $y(n)=-x(n)$ for each $n$, as desired. 
		
		The main step in the proof is verifying the following statement holds for each $n$: 
		\begin{enumerate}
			\item For each $F\subset \nn$ with $|F|=n$ and $(\varepsilon_i)_{i\in F} \subset \{\pm 1\}$, $\sum_{i \in F} \varepsilon_i V(e_i)\in S_X$.
			\item For each $F\subset \nn$ with $|F|=n$ and  $(a_i)_{i\in F}$ with $\sum_{i \in F} a_i V(e_i)\in S_X$ there is an $x \in c_0$ with $\supp\, x\cap F=\emptyset $ so that $V(\sum_{i\in F} a_i e_i +x)= \sum_{i\in F} a_i V(e_i)$.
		\end{enumerate} 
		
		We proceed by induction on $n$. The base step $n=1$ follows from the definitions and the fact that $V(-x)=-V(x)$. Assume the statement holds for each $k<n$, fix $F=\{m_1, \ldots, m_n\}$ and $(\varepsilon_i) \subset\{\pm 1\}$. By the induction hypothesis there exist vectors $x,y\in c_0$ with $\{m_1, \ldots, m_{n-1}\}\cap \supp\, x=\emptyset$ and $\{m_2, \ldots, m_n\}\cap \supp\, y=\emptyset$ satisfying
		$$V(\sum_{i=1}^{n-1}\varepsilon_{m_i} e_{m_i} +x)=\sum_{i=1}^{n-1}\varepsilon_{m_i} Ve_{m_i}$$
		and 
		$$V(-\sum_{i=2}^{n-1}\varepsilon_{m_i} e_{m_i} + \varepsilon_{m_n} e_{m_n} +y)=-\sum_{i=2}^{n-1}\varepsilon_{m_i} V(e_{m_i}) + \varepsilon_{m_n} V(e_{m_n}).$$
		Then 
		\begin{equation}
		\begin{split}
		1 & =\|\varepsilon_{m_1}V(e_{m_1})+\varepsilon_{m_n}V(e_{m_n})\| \\
		&  =\|V(\sum_{i=1}^{n-1}\varepsilon_{m_i} e_{m_i} +x) + V(-\sum_{i=2}^{n-1}\varepsilon_{m_i} e_{m_i} + \varepsilon_{m_n} e_{m_n} +y)\|  \\
		& = \| \varepsilon_{m_1} e_{m_1} + \varepsilon_{m_n} e_{m_n} +x +y\| \geqslant |\varepsilon_{m_n} + x(m_n) + y(m_n)|
		\end{split}
		\end{equation}
		Since $y(m_n)=0$ we have  $|\varepsilon_{m_n} + x(m_n)|\leqslant 1$. Consequently,
		\begin{equation}
		\begin{split}
		\|\sum_{i=1}^{n}\varepsilon_{m_i} Ve_{m_i}\| & = \|V(\sum_{i=1}^{n-1}\varepsilon_{m_i} e_{m_i} +x) + \varepsilon_{m_n} Ve_{m_n}\|\\
		& = \|\sum_{i=1}^{n}\varepsilon_{m_i} e_{m_i} +x\| =\max\{1,|\varepsilon_{m_n} + x(m_n)|,\|x\|\}= 1.   
		\end{split}
		\end{equation}
		This proves item (1) of the induction for $n$. For the second part we fix $(a_i)$ with $\sum_{1=1}^n a_{m_i}V(e_{m_i})\in S_X$. Let $\sum_i b_i e_i\in S_{c_0}$ so that $V(\sum_i b_i e_i)=\sum_{1=1}^n a_{m_i}V(e_{m_i})$. Let $x=\sum_{i\not\in F} b_ie_i$. Then, trivially, 
		$$V(\sum_{i=1}^n b_{m_i} e_{m_i}+x)=\sum_{1=1}^n a_{m_i}V(e_{m_i}).$$
		It suffices to show the $b_{m_i}=a_{m_i}$ for each $1\leqslant i \leqslant n$. Fix $j\in \{1,\ldots,n\}$. Then 
		by Lemma \ref{allscalars} we have
		$$\| \sum_{i=1}^na_{m_i} V(e_{m_i})+\sgn(b_{m_j})e_{m_j}\| = |a_{m_j}+\sgn(b_{m_j})|$$
		and 
		$$\| \sum_{i=1}^nb_i e_i+\sgn(b_{m_j})e_{m_j}\| = |b_{m_j}+\sgn(b_{m_j})|=1+|b_{m_j}|.$$
		As the above are equal, we see that $a_{m_j}=b_{m_j}$. The induction is complete.
		
		We claim that the vector $x$ in item (2) is 0 for any $\sum_{i \in F} a_i V(e_i)\in S_X$. Let $F\subset \mathbb{N}$ and find $x$ so that,  $$V(\sum_{i \in F} a_i e_i+x)=\sum_{i \in F} a_i V(e_i)$$
		Let $j\not\in F$. Then 
		$$|x(j)+\sgn(x(j))|\leqslant \|\sum_{i \in F} a_i e_i+x + \sgn(x(j))e_j\|=\|\sum_{i \in F} a_i V(e_i) + \sgn(x(j))V(e_j)\|= 1$$
		Then $x(j)=0$, and so $x=0$.
		
		Therefore, for $\sum_i a_i e_i\in S_{c_0}$ we have proved $V(\sum_i a_i e_i)= \sum_i a_i V(e_i)$. Let $U:c_0 \to X$ be defined by $U(\sum_i b_i e_i) = \sum_i b_i V(e_i).$ 
		By Lemma \ref{allscalars} we that $U$ is an isometry on a dense subset of $c_0$ and thus extends to an isometry on all of $c_0$. The results follow. \end{proof}

	\bibliographystyle{abbrv}
	\bibliography{bib_source}
	
\end{document}